\documentclass[11pt]{amsart}

\usepackage{graphicx}
\usepackage{amssymb}
\usepackage{amsfonts}
\usepackage[]{amsmath}
\usepackage[]{epsfig}
\usepackage[]{pstricks}
\usepackage[]{float}
\newpsobject{malla}{psgrid}{subgriddiv=1,griddots=10,gridlabels=6pt}

\newtheorem{theorem}{Theorem}[section]
\newtheorem{lemma}{Lemma}[section]

\newtheorem{remark}{Remark}[section]
\newtheorem{proposition}{Proposition}[section]

\numberwithin{equation}{section}

\newcommand{\T}{\mathbb T}

\newcommand{\ov}{\overline}

\newcommand{\la}{\langle}
\newcommand{\ra}{\rangle}

\newcommand{\p}{\partial}

\begin{document}

\title[Global solutions for the periodic NLS-KdV system]{Global well-posedness
for a NLS-KdV system on $\mathbb{T}$}

\author{C. Matheus}
\address{Carlos Matheus} 
\address{IMPA, Estrada Dona Castorina 110, Rio de Janeiro, 22460--320,
Brazil.} \email{matheus@impa.br}


\date{November 16, 2005}

\begin{abstract}We prove that the Cauchy problem of the Schr\"odinger - Korteweg
- deVries (NLS-KdV) system on $\mathbb{T}$ is globally well-posed for initial 
data $(u_0,v_0)$ below the energy space $H^1\times H^1$. More precisely, we 
show that the non-resonant NLS-KdV is globally well-posed for initial data
$(u_0,v_0)\in H^s(\mathbb{T})\times H^s(\mathbb{T})$ with $s>11/13$ and the
resonant NLS-KdV is globally well-posed for initial data 
$(u_0,v_0)\in H^s(\mathbb{T})\times H^s(\mathbb{T})$ with $s>8/9$. The idea of
the proof of this theorem is to apply the I-method of Colliander, Keel,
Staffilani, Takaoka and Tao in order to improve the results of Arbieto, Corcho
and Matheus concerning the global well-posedness of the NLS-KdV on $\mathbb{T}$
in the energy space $H^1\times H^1$.
\end{abstract}

\maketitle


\section{Introduction}\label{s.intro}
We consider the Cauchy problem of the Schr\"odinger-Korteweg-deVries system  
\begin{equation}\label{e.nls-kdv}
\begin{cases}
i\partial_tu + \partial_x^2u = \alpha uv + \beta |u|^2u,\\
\partial_tv + \partial_x^3v + \tfrac{1}{2}\partial_x(v^2) = 
\gamma \partial_x(|u|^2),\\
u(x,0)=u_0(x),\;v(x,0)=v_0(x), &  t\in\mathbb{R}.
\end{cases} 
\end{equation}

The Schr\"odinger-Korteweg-deVries (NLS-KdV) system  naturally appears in fluid 
mechanics and plasma physics as a model of interaction between a short-wave 
$u=u(x,t)$ and a long-wave $v=v(x,t)$. 

In this paper we are interested in global solutions of the NLS-KdV system for
rough initial data. Before stating our main results, let us recall some of the
recent theorems of local and global well-posedness theory of the Cauchy
problem~(\ref{e.nls-kdv}). 

For continuous spatial variable (i.e., $x\in\mathbb{R}$), Corcho and
Linares~\cite{Corcho}
recently proved that the NLS-KdV system is locally well-posed for initial data
$(u_0,v_0)\in H^k(\mathbb{R})\times H^s(\mathbb{R})$ with $k\geq 0$, $s>-3/4$
and 
\begin{itemize}
\item $k-1\leq s\leq 2k-1/2$ if $k\leq 1/2$ , 
\item $k-1\leq s<k+1/2$ if $k>1/2$.
\end{itemize}
Furthermore, they were able to prove the global well-posedness of the NLS-KdV 
system in the energy $H^1\times H^1$ using three conserved quantities 
discovered by M. Tsutsumi~\cite{MTsutsumi}, whenever $\alpha\gamma>0$. 

Also, Pecher~\cite{Pecher} improved this global well-posedness
result by an application of the I-method of Colliander, Keel, Stafillani,
Takaoka and Tao (for instance, see~\cite{CKSTT1}) combined with some refined 
bilinear estimates. In particular, Pecher proved that, if $\alpha\gamma>0$, 
the NLS-KdV system is globally well-posed for initial data 
$(u_0,v_0)\in H^s\times H^s$ with $s>3/5$ in the resonant case $\beta=0$ and 
$s>2/3$ in the non-resonant case $\beta\neq 0$.

On the other hand, in the periodic setting (i.e., $x\in\mathbb{T}$), Arbieto,
Corcho and Matheus~\cite{ACM} proved the local well-posedness of the NLS-KdV
system for initial data $(u_0,v_0)\in H^k\times H^s$ with $0\leq s\leq 4k-1$ and
$-1/2\leq k-s\leq 3/2$. Also, using the same three conserved quantities
discovered by M. Tsutsumi, one obtains the global well-posedness of NLS-KdV on
$\mathbb{T}$ in the energy space $H^1\times H^1$ whenever $\alpha\gamma>0$.

Motivated by this scenario, we combine the new bilinear estimates of Arbieto, 
Corcho and Matheus~\cite{ACM} with the I-method of Tao and his collaborators to
prove the following result 

\begin{theorem}\label{t.A}The NLS-KdV system~(\ref{e.nls-kdv}) on $\mathbb{T}$ 
is globally well-posed for initial data $(u_0,v_0)\in H^s(\mathbb{T})\times
H^s(\mathbb{T})$ with $s>11/13$ in the non-resonant case $\beta\neq 0$ and
$s>8/9$ in the resonant case $\beta=0$, whenever $\alpha\gamma>0$.
\end{theorem}

The paper is organized as follows. In the section~\ref{s.preliminaries}, we
discuss the preliminaries for the proof of the theorem~\ref{t.A}:
Bourgain spaces and its properties, linear estimates, standard estimates for 
the non-linear terms $|u|^2 u$ and $\p_x (v^2)$, the bilinear estimates of
Arbieto, Corcho and Matheus~\cite{ACM} for the coupling terms $uv$ and
$\p_x(|u|^2)$, the I-operator and its properties. In the section~\ref{s.local},
we apply the results of the section~\ref{s.preliminaries} to get a variant of
the local well-posedness result of~\cite{ACM}. In the
section~\ref{s.conservation}, we recall some conserved quantities
of~(\ref{e.nls-kdv}) and its modification by the introduction of the I-operator;
moreover, we prove that two of these modified energies are almost 
conserved. Finally, in the section~\ref{s.global}, we combine the almost
conservation results in section~\ref{s.conservation} with the local
well-posedness result in section~\ref{s.local} to conclude the proof of the
theorem~\ref{t.A}.  

\section{Preliminaries}\label{s.preliminaries}

A sucessful procedure to solve some dispersive equations (such as the
nonlinear Schr\"odinger and KdV equations) is to use the Picard's fixed point 
method in the following spaces: 
 
\begin{equation*}
\begin{split}
\|f\|_{X^{k,b}}&:= \left(\int\sum\limits_{n\in\mathbb{Z}} \la n\ra^{2k}
\la\tau+n^2\ra^{2b}|\widehat{f}(n,\tau)| d\tau\right)^{1/2} \\ 
&= \|U(-t) f\|_{H_t^b(\mathbb{R},H_x^k)}
\end{split}
\end{equation*}

\begin{equation*}
\begin{split}
\|g\|_{Y^{s,b}}&:= \left(\int\sum\limits_{n\in\mathbb{Z}} \la n\ra^{2s}
\la\tau-n^3\ra^{2b}|\widehat{g}(n,\tau)| d\tau\right)^{1/2} \\ 
&= \|V(-t) f\|_{H_t^b(\mathbb{R},H_x^s)}
\end{split}
\end{equation*}
where $\la\cdot\ra:= 1+ |\cdot|$, $U(t)=e^{it\p_x^2}$ and $V(t)=e^{-t\p_x^3}$. 
These spaces are called Bourgain spaces. Also, 
we introduce the restriction in time norms  

\begin{equation*}
\|f\|_{X^{k,b}(I)}:=\inf\limits_{\widetilde{f}|_I=f} \|\widetilde{f}\|_{X^{k,b}}
\quad \textrm{ and } \quad 
\|g\|_{Y^{s,b}(I)}:=\inf\limits_{\widetilde{g}|_I=g} \|\widetilde{g}\|_{Y^{s,b}}
\end{equation*}
where $I$ is a time interval.

The interaction of the Picard method has been based around the spaces 
$Y^{s,1/2}$. Because we are interested in the continuity of the flow associated
to~(\ref{e.nls-kdv}) and the $Y^{s,1/2}$ norm do not control the
$L_t^{\infty}H_x^s$ norm, we modify the Bourgain spaces as follows:

\begin{equation*}
\|u\|_{X^k}:= \|u\|_{X^{k,1/2}} + \|\la n\ra^k\widehat{u}(n,\tau)\|_{L_n^2
L_{\tau}^1}
\quad \textrm{ and } \quad 
\|v\|_{Y^{s}}:=\|v\|_{Y^{s,1/2}} + \|\la n\ra^s\widehat{v}(n,\tau)\|_{L_n^2
L_{\tau}^1}
\end{equation*}
and, given a time interval $I$, we consider the restriction in time of the
$X^k$ and $Y^s$ norms 
\begin{equation*}
\|u\|_{X^{k}(I)}:=\inf\limits_{\widetilde{u}|_I=u} \|\widetilde{u}\|_{X^{k}}
\quad \textrm{ and } \quad 
\|v\|_{Y^{s}(I)}:=\inf\limits_{\widetilde{v}|_I=v} \|\widetilde{v}\|_{Y^{s}}
\end{equation*}

Furthermore, the mapping properties of $U(t)$ and $V(t)$ naturally leads one to
consider the companion spaces
\begin{equation*}
\|u\|_{Z^k}:=\|u\|_{X^{k,-1/2}}+
\left\|\frac{\la n\ra^k\widehat{u}(n,\tau)}{\la\tau+n^2\ra}\right\|_{L_n^2
L_{\tau}^1} \quad \textrm{ and } \quad 
\|v\|_{W^s}:=\|v\|_{Y^{s,-1/2}}+
\left\|\frac{\la n\ra^s\widehat{v}(n,\tau)}{\la\tau-n^3\ra}\right\|_{L_n^2
L_{\tau}^1}
\end{equation*} 

In the sequel, $\psi$ denotes a non-negative smooth bump function supported on
$[-2,2]$ with $\psi=1$ on $[-1,1]$ and $\psi_{\delta}(t):=\psi(t/\delta)$ for
any $\delta>0$.

Next, we recall some properties of the Bourgain spaces:

\begin{lemma}\label{l.Strichartz} $X^{0,3/8}([0,1]), Y^{0,1/3}([0,1])\subset 
L^4(\mathbb{T}\times [0,1])$. More precisely, 
\begin{equation*}
\|\psi(t) f\|_{L_{xt}^4}\lesssim \|f\|_{X^{0,3/8}} \quad \textrm{and} \quad
\|\psi(t) g\|_{L_{xt}^4}\lesssim \|g\|_{Y^{0,1/3}}.
\end{equation*} 
\end{lemma}

\begin{proof}See~\cite{Bourgain}.
\end{proof}

Another basic property of these spaces are their stability under time
localization:

\begin{lemma}\label{l.time-localization} Let $X^{s,b}_{\tau = h(\xi)}:= \{f :
\la\tau-h(\xi)\ra^b\la\xi\ra^s |\widehat{f}(\tau,\xi)|\in L^2\}$. Then, 
\begin{equation*}
\|\psi(t) f\|_{X^{s,b}_{\tau=h(\xi)}}\lesssim_{\psi,b}
\|f\|_{X^{s,b}_{\tau=h(\xi)}} 
\end{equation*}
for any $s,b\in\mathbb{R}$. Moreover, if $-1/2<b'\leq b<1/2$, then for any
$0<T<1$, we have  
\begin{equation*}
\|\psi_T(t) f\|_{X_{\tau=h(\xi)}^{s,b'}}\lesssim_{\psi,b',b} T^{b-b'}
\|f\|_{X^{s,b}_{\tau=h(\xi)}}.
\end{equation*}
\end{lemma}

\begin{proof}First of all, note that $\la\tau-\tau_0-h(\xi)\ra^{b}\lesssim_b 
\la\tau_0\ra^{|b|}\la\tau-h(\xi)\ra^{b}$, 
from which we obtain 
$$\|e^{it\tau_0}f\|_{X_{\tau=h(\xi)}^{s,b}}\lesssim_b \la\tau_0\ra^{|b|} 
\|f\|_{X_{\tau=h(\xi)}^{s,b}}.$$  
Using that $\psi(t)=\int\widehat{\psi}(\tau_0) e^{it\tau_0}d\tau_0$, we 
conclude 
$$\|\psi(t)f\|_{X_{\tau=h(\xi)}^{s,b}}\lesssim_b 
\left(\int|\widehat{\psi}(\tau_0)| \la\tau_0\ra^{|b|}\right) 
\|f\|_{X_{\tau=h(\xi)}^{s,b}}.$$  
Since $\psi$ is smooth with compact support, the first estimate follows.  

Next we prove the second estimate. By conjugation we may assume $s=0$ and, 
by composition it suffices to treat the cases $0\leq b'\leq b$ or 
$\leq b'\leq b\leq 0$. By duality, we may take $0\leq b'\leq b$. 
Finally, by interpolation with the trivial case $b'=b$, we may consider 
$b'=0$. This reduces matters to show that 
$$\|\psi_T(t)f\|_{L^2}\lesssim_{\psi,b} T^b\|f\|_{X_{\tau=h(\xi)}^{0,b}}$$ 
for $0<b<1/2$. Partitioning the frequency spaces into the cases 
$\la\tau-h(\xi)\ra\geq 1/T$ and $\la\tau-h(\xi)\leq 1/T$, we see that in the 
former case we'll have 
$$\|f\|_{X_{\tau=h(\xi)}^{0,0}}\leq T^b\|f\|_{X_{\tau=h(\xi)}^{0,b}}$$ 
and the desired estimate follows because the multiplication by $\psi$ is a 
bounded operation in Bourgain's spaces. In the latter case, by Plancherel and 
Cauchy-Schwarz 
\begin{equation*} 
\begin{split} 
\|f(t)\|_{L_x^2}&\lesssim \|\widehat{f(t)}(\xi)\|_{L_{\xi}^2} \lesssim 
\left\|\int_{\la\tau-h(\xi)\ra\leq 1/T}|\widehat{f}(\tau,\xi)|d\tau) 
\right\|_{L_{\xi}^2} \\ &\lesssim_b T^{b-1/2} 
\left\|\int\la\tau-h(\xi)\ra^{2b} |\widehat{f}(\tau,\xi)|^2 
d\tau)^{1/2}\right\|_{L_{\xi}^2} = T^{b-1/2}\|f\|_{X_{\tau=h(\xi)}^{s,b}}. 
\end{split} 
\end{equation*}  
Integrating this against $\psi_T$ concludes the proof of the lemma. 
\end{proof}

Also, we have the following duality relationship between $X^k$ (resp., $Y^s$) 
and $Z^k$ (resp., $W^s$):
  
\begin{lemma}\label{l.duality} We have 
$$\left|\int \chi_{[0,1]}(t) f(x,t) g(x,t) dx dt\right|\lesssim \|f\|_{X^s}
\|g\|_{Z^{-s}}$$
and
$$\left|\int \chi_{[0,1]}(t) f(x,t) g(x,t) dx dt\right|\lesssim \|f\|_{Y^s}
\|g\|_{W^{-s}}$$
for any $s$ and any $f, g$ on $\mathbb{T}\times\mathbb{R}$.
\end{lemma}

\begin{proof}See~\cite[p. 182--183]{CKSTT} (note that, although this result is
stated only for the spaces $Y^s$ and $W^s$, the same proof adapts for the spaces
$X^k$ and $Z^k$).
\end{proof}

Now, we recall some linear estimates related to the semigroups $U(t)$ and
$V(t)$:

\begin{lemma}[Linear estimates]\label{l.linear}It holds 
\begin{itemize}
\item $\|\psi(t)U(t)u_0\|_{Z^k}\lesssim \|u_0\|_{H^k}$ and 
$\|\psi(t)V(t)v_0\|_{W^s}\lesssim \|v_0\|_{H^s}$;
\item $\|\psi_{T}(t)\int_0^t U(t-t') F(t') dt'\|_{X^k}\lesssim \|F\|_{Z^k}$ and 
$\|\psi_{T}(t)\int_0^t V(t-t') G(t') dt'\|_{Y^s}\lesssim \|G\|_{W^s}$.
\end{itemize}
\end{lemma}

\begin{proof}See~\cite{CKSTT1},~\cite{CKSTT} or~\cite{ACM}.
\end{proof}

Furthermore, we have the following well-known multiinear estimates for the cubic
term $|u|^2 u$ of the nonlinear Schr\"odinger equation and the nonlinear term
$\p_x(v^2)$ of the KdV equation:

\begin{lemma}\label{l.u2u}
$\|uv\ov{w}\|_{Z^k}\lesssim \|u\|_{X^{k,\frac{3}{8}}}\|v\|_{X^{k,\frac{3}{8}}}
\|w\|_{X^{k,\frac{3}{8}}}$ for any $k\geq 0$.
\end{lemma}

\begin{proof}See~\cite{Bourgain} and~\cite{ACM}.
\end{proof}

\begin{lemma}\label{l.dv2}
$\|\p_x(v_1 v_2)\|_{W^s}\lesssim
\|v_1\|_{Y^{s,\frac{1}{3}}}\|v_2\|_{Y^{s,\frac{1}{2}}}+
\|v_1\|_{Y^{s,\frac{1}{2}}}\|v_2\|_{Y^{s,\frac{1}{3}}}$ for any $s\geq -1/2$, if
$v_1=v_1(x,t)$ and $v_2=v_2(x,t)$ are $x$-periodic functions having zero
$x$-mean for all $t$.
\end{lemma}

\begin{proof}See~\cite{Bourgain},~\cite{CKSTT1} and~\cite{ACM}.
\end{proof}

Next, we revisit the bilinear estimates of mixed Schr\"odinger-Airy type of
Arbieto, Corcho and Matheus~\cite{ACM} for the coupling terms $uv$ and $\p_x(|u|^2)$ of the
NLS-KdV system.

\begin{lemma}\label{l.uv}
$\|uv\|_{Z^k}\lesssim \|u\|_{X^{k,\frac{3}{8}}}\|v\|_{Y^{s,\frac{1}{2}}} + 
\|u\|_{X^{k,\frac{1}{2}}}\|v\|_{Y^{s,\frac{1}{3}}}$ whenever $s\geq 0$ and 
$k-s\leq 3/2$.
\end{lemma}

\begin{lemma}\label{l.du2}
$\|\p_x(u_1\ov{u_2})\|_{W^s}\lesssim \|u_1\|_{X^{k,3/8}}\|u_2\|_{X^{k,1/2}} + 
\|u_1\|_{X^{k,1/2}}\|u_2\|_{X^{k,3/8}}$ whenever $1+s\leq 4k$ and $k-s\geq
-1/2$.
\end{lemma}

\begin{remark}Although the lemmas~\ref{l.uv} and~\ref{l.du2} are not stated as
above in~\cite{ACM}, it is not hard to obtain them from the 
calculations of Arbieto, Corcho and Matheus.  
\end{remark}

Finally, we introduce the I-operator: let $m(\xi)$ be a smooth non-negative
symbol on $\mathbb{R}$ which equals $1$ for $|\xi|\leq 1$ and equals
$|\xi|^{-1}$ for $|\xi|\geq 2$. For any $N\geq 1$ and $\alpha\in\mathbb{R}$,
denote by $I_N^{\alpha}$ the spatial Fourier multiplier 
\begin{equation*}
\widehat{I_N^{\alpha}f}(\xi) = m\left(\frac{\xi}{N}\right)^{\alpha}
\widehat{f}(\xi).
\end{equation*}

For latter use, we recall the following general interpolation lemma:
\begin{lemma}[Lemma 12.1 of~\cite{CKSTT}]\label{l.interpolation}Let $\alpha_0>0$
and $n\geq 1$. Suppose $Z, X_1, X_2,\dots, X_n$ are translation-invariant Banach
spaces and $T$ is a translation invariant $n$-linear operator such that 
\begin{equation*}
\|I_1^{\alpha} T(u_1,\dots, u_n)\|_{Z}\lesssim \prod\limits_{j=1}^n
\|I_1^{\alpha}u_j\|_{X_j},
\end{equation*}
for all $u_1,\dots,u_n$ and $0\leq\alpha\leq\alpha_0$. Then, 
\begin{equation*}
\|I_N^{\alpha} T(u_1,\dots, u_n)\|_{Z}\lesssim \prod\limits_{j=1}^n
\|I_N^{\alpha}u_j\|_{X_j}
\end{equation*}
for all $u_1,\dots,u_n$, $0\leq\alpha\leq\alpha_0$ and $N\geq 1$. Here the
implicit constant is independent of $N$.
\end{lemma}

After these preliminaries, we can proceed to the next section where 
a variant of the local well-posedness of Arbieto, Corcho and 
Matheus is obtained. 

In the sequel we take $N\gg 1$ a large integer and denote by $I$ the operator
$I=I_N^{1-s}$ for a given $s\in\mathbb{R}$.

\section{A variant local well-posedness result}\label{s.local} 

This section is devoted to the proof of the following proposition:

\begin{proposition}\label{p.local}For any $(u_0,v_0)\in H^s(\mathbb{T})\times 
H^s(\mathbb{T})$ with $\int_{\mathbb{T}} v_0 = 0$ and $s\geq 1/3$, the periodic 
NLS-KdV system~(\ref{e.nls-kdv})
has a unique local-in-time solution on the time interval $[0,\delta]$ for some
$\delta\leq 1$ and 
\begin{equation}\label{e.local}
\delta\sim 
\begin{cases} 
(\|Iu_0\|_{X^1}+\|Iv_0\|_{Y^1})^{-\frac{16}{3}-}, \textrm{ if } \beta\neq 0, \\ 
(\|Iu_0\|_{X^1}+\|Iv_0\|_{Y^1})^{-8-}, \textrm{ if } \beta = 0.
\end{cases}
\end{equation}
Moreover, we have $\|Iu\|_{X^1}+\|Iv\|_{Y^1}\lesssim
\|Iu_0\|_{X^1}+\|Iv_0\|_{Y^1}$.
\end{proposition}

\begin{proof}We apply the I-operator to the NLS-KdV system~(\ref{e.nls-kdv}) so
that 
\begin{equation*}
\begin{cases}
i Iu_t + Iu_{xx} = \alpha I(uv) + \beta I(|u|^2 u), \\ 
Iv_t + Iv_{xxx} + I(v v_x) = \gamma I(|u|^2)_x, \\ 
Iu(0) = Iu_0, \ Iv(0) = Iv_0.
\end{cases}
\end{equation*}
To solve this equation, we seek for some fixed point of the integral maps 
\begin{equation*}
\begin{split}
\Phi_1(Iu, Iv):= U(t) Iu_0 -i\int_0^t U(t-t') \{\alpha I(u(t')v(t')) + \beta
I(|u(t')|^2 u(t'))\} dt', \\ 
\Phi_2(Iu, Iv):= V(t)Iv_0 -\int_0^t V(t-t')\{I(v(t')v_x(t')) - \gamma
I(|u(t')|^2)_x\} dt'.
\end{split}
\end{equation*}

The interpolation lemma~\ref{l.interpolation} applied to the linear and
multilinear estimates in the lemmas~\ref{l.linear},~\ref{l.u2u},~\ref{l.dv2},
~\ref{l.uv} and~\ref{l.du2} yields, in view of the
lemma~\ref{l.time-localization},  
\begin{equation*}
\begin{split}
\|\Phi_1(Iu, Iv)\|_{X^1}\lesssim \|Iu_0\|_{H^1} + \alpha\delta^{\frac{1}{8}-}
\|Iu\|_{X^1}\|Iv\|_{Y^1} + \beta\delta^{\frac{3}{8}-}\|Iu\|_{X^1}^3, \\
\|\Phi_2(Iu, Iv)\|_{Y^1}\lesssim \|Iv_0\|_{H^1} + 
\delta^{\frac{1}{6}-}\|Iv\|_{Y^1}^2 + 
\gamma\delta^{\frac{1}{8}-}\|Iu\|_{X^1}^2.
\end{split}
\end{equation*}
In particular, these integrals maps are contractions provided that
$\beta\delta^{\frac{3}{8}-}(\|Iu_0\|_{H^1}+\|Iv_0\|_{H^1})^2 \ll 1$ and
$\delta^{\frac{1}{8}-}(\|Iu_0\|_{H^1}+\|Iv_0\|_{H^1})\ll 1$. This completes the
proof of the proposition~\ref{p.local}.
\end{proof}

\section{Modified energies}\label{s.conservation} 

Consider the following three quantities:

\begin{equation}\label{e.mass}
M(u):=\|u\|_{L^2},
\end{equation}

\begin{equation}\label{e.momentum}
L(u,v):= \alpha \|v\|_{L^2}^2 + 2\gamma\int \Im (u\ov{u_x}) dx,
\end{equation}

\begin{equation}\label{e.energy}
E(u,v):= \alpha\gamma\int v |u|^2 dx + \gamma\|u_x\|_{L^2}^2 + \frac{\alpha}{2}
\|v_x\|_{L^2}^2 -\frac{\alpha}{6}\int v^3 dx + \frac{\beta\gamma}{2}\int |u|^4
dx.
\end{equation}

In the sequel, we suppose $\alpha\gamma>0$. Note that 

\begin{equation}\label{e.L1}
|L(u,v)|\lesssim \|v\|_{L^2}^2 + M \|u_x\|_{L^2}
\end{equation} 
and 
\begin{equation}\label{e.L2}
\|v\|_{L^2}^2\lesssim |L| + M \|u_x\|_{L^2}.
\end{equation}

Also, the Gagliardo-Nirenberg and Young inequalities implies 

\begin{equation}\label{e.E1}
\|u_x\|_{L^2}^2+\|v_x\|_{L^2}^2\lesssim |E| + |L|^{\frac{5}{3}} + M^8 + 1
\end{equation} 
and 
\begin{equation}\label{e.E2}
|E|\lesssim \|u_x\|_{L^2}^2 + \|v_x\|_{L^2}^2 + |L|^{\frac{5}{3}} + M^8 + 1
\end{equation}

In particular, combining the bounds~(\ref{e.L1}) and~(\ref{e.E2}), 

\begin{equation}\label{e.E3}
|E|\lesssim \|u_x\|_{L^2}^2 + \|v_x\|_{L^2}^2 + \|v\|_{L^2}^{\frac{10}{3}} +
M^{10} + 1.
\end{equation}

Moreover, from the bounds~(\ref{e.L2}) and~(\ref{e.E1}), 

\begin{equation}\label{e.E4}
\|v\|_{L^2}^2\lesssim |L| + M |E|^{1/2} + M^6 + 1
\end{equation}
and hence 
\begin{equation}\label{e.E5}
\|u\|_{H^1}^2 + \|v\|_{H^1}^2\lesssim |E| + |L|^{5/3} + M^8 + 1
\end{equation}

\begin{equation}\label{e.al}
\begin{split}
\frac{d}{dt} L(Iu, Iv) &= 2\alpha\int Iv (Iv Iv_x - I(v v_x)) dx +
2\alpha\gamma\int Iv (I(|u|^2)-|Iu|^2)_x dx \\ 
&+ 4\alpha\gamma\Re\int I\ov{u}_x (Iu Iv - I(uv)) dx  
+ 4\beta\gamma\Re\int ((Iu)^2 I\ov{u} - I(u^2\ov{u})) I\ov{u}_x dx \\ 
&=: \sum\limits_{j=1}^4 L_j.
\end{split}
\end{equation}

\begin{equation}\label{e.ae}
\begin{split}
\frac{d}{dt} E(Iu, Iv) &= \alpha\int (I(vv_x)-IvIv_x)Iv_{xx} dx + 
\frac{\alpha}{2}\int (Iv)^2 (I(vv_x)-IvIv_x) dx + \\ 
&+ 2\beta\gamma\Im\int (I(|u|^2 u)_x - ((Iu)^2 I\ov{u})_x) I\ov{u}_x dx \\ 
&+ \alpha\gamma\int |Iu|^2 (Iv Iv_x - I(v v_x)) dx + \alpha\gamma\int (|Iu|^2 -
I(|u|^2))Iv Iv_x dx \\ 
&+ \alpha\gamma\int Iv_{xx} (|Iu|^2-I(|u|^2))_x dx -2\alpha\gamma\Im\int Iu_x
(I(\ov{u}v)-I\ov{u} Iv)_x dx \\ 
&+ \alpha\gamma^2 \int (I(|u|^2) - |Iu|^2)_x |Iu|^2 dx + 2\alpha^2\gamma\Im\int
Iv Iu (I(\ov{u}v - I\ov{u} Iv)) dx\\ 
&+ 2\beta^2\gamma\Im\int Iu(I\ov{u})^2 (I(|u|^2 u) - (Iu)^2 I\ov{u}) dx \\ 
&-2\alpha\beta\gamma\Im\int Iv Iu (I(|u|^2 \ov{u}) - Iu (I\ov{u}))^2 dx
-2\alpha\beta\gamma\Im\int (Iu)^2 I\ov{u} (I(\ov{u}v) - I\ov{u} Iv) dx \\ 
&=: \sum\limits_{j=1}^{12} E_j
\end{split}
\end{equation}

\subsection{Estimates for the modified L-functional}\label{s.l}

\begin{proposition}\label{p.al}Let $(u,v)$ be a solution of~(\ref{e.nls-kdv}) on
the time interval $[0,\delta]$. Then, for any $N\geq 1$ and $s>1/2$, 
\begin{equation}\label{e.aL}
\begin{split}
&|L(Iu(\delta), Iv(\delta)) - L(Iu(0), Iv(0))|\lesssim \\ 
&N^{-1+}\delta^{\frac{19}{24}-} (\|Iu\|_{X^{1,1/2}}+\|Iv\|_{Y^{1,1/2}})^3 + 
N^{-2+}\delta^{\frac{1}{2}-}\|Iu\|_{X^{1,1/2}}^4.
\end{split}
\end{equation}
\end{proposition}

\begin{proof}Integrating~(\ref{e.al}) with respect to $t\in [0,\delta]$, it
follows that we have to bound the (integral over $[0,\delta]$ of the) four 
terms on the right hand side. To simplify the computations, we assume that the
Fourier transform of the functions are non-negative and we ignore the appearance
of complex conjugates (since they are irrelevant in our subsequent arguments).
Also, we make a dyadic decomposition of the frequencies $|n_i|\sim N_j$ in many
places. In particular, it will be important to get extra factors $N_j^{0-}$
everywhere in order to sum the dyadic blocks. 

We begin with the estimate of $\int_{0}^{\delta} L_1$. It is sufficient to show
that 
\begin{equation}\label{e.aL1}
\begin{split}
&\int_{0}^{\delta}\sum\limits_{n_1+n_2+n_3=0}
\Big|\frac{m(n_1+n_2) - m(n_1) m(n_2)}{m(n_1) m(n_2)}\Big|
\widehat{v_1}(n_1,t) |n_2|\widehat{v_2}(n_2,t) \widehat{v_3}(n_3,t) 
\lesssim \\ 
&N^{-1}
\delta^{\frac{5}{6}-}\prod\limits_{j=1}^3 \|v_j\|_{Y^{1,1/2}}
\end{split}
\end{equation}
\begin{itemize}
\item $|n_1|\ll |n_2|\sim |n_3|$, $|n_2|\gtrsim N$. In this case, note that 
\begin{equation*}
\begin{cases}
\big|\frac{m(n_1+n_2) - m(n_1) m(n_2)}{m(n_1) m(n_2)}\big|\lesssim 
|\frac{\nabla m(n_2)\cdot n_1}{m(n_2)}|\lesssim \frac{N_1}{N_2}, 
\textrm{ if } |n_1|\leq N, \text{ and }\\
\big|\frac{m(n_1+n_2) - m(n_1) m(n_2)}{m(n_1) m(n_2)}\big|\lesssim 
\left(\frac{N_1}{N}\right)^{1/2}, \textrm{ if } |n_1|\geq N.
\end{cases} 
\end{equation*}
Hence, using the lemmas~\ref{l.Strichartz} and~\ref{l.time-localization} we
obtain 
\begin{equation*}
|\int_{0}^{\delta} L_1|\lesssim \frac{N_1}{N_2} \|v_1\|_{L^4} \|(v_2)_x\|_{L^4}
\|v_3\|_{L^2}\lesssim N^{-2+}\delta^{\frac{5}{6}-} N_{\max}^{0-}
\prod\limits_{j=1}^3\|v_i\|_{Y^{1,1/2}}
\end{equation*} 
if $|n_1|\leq N$, and 
\begin{equation*}
|\int_{0}^{\delta} L_1|\lesssim \left(\frac{N_1}{N}\right)^{1/2} \frac{1}{N_1
N_3} \delta^{\frac{5}{6}-}\prod\limits_{j=1}^3\|v_i\|_{Y^{1,1/2}}\lesssim
N^{-2+}\delta^{\frac{5}{6}-} N_{\max}^{0-}
\prod\limits_{j=1}^3\|v_i\|_{Y^{1,1/2}}.
\end{equation*}
\item $|n_2|\ll |n_1|\sim |n_3|$, $|n_1|\gtrsim N$. This case is similar to the
previous one. 
\item $|n_1|\sim |n_2|\gtrsim N$. The multiplier is bounded by 
\begin{equation*}
\big|\frac{m(n_1+n_2) - m(n_1) m(n_2)}{m(n_1) m(n_2)}\big|\lesssim 
\left(\frac{N_1}{N}\right)^{1-}.
\end{equation*}
In particular, using the lemmas~\ref{l.Strichartz}
and~\ref{l.time-localization}, 
\begin{equation*}
|\int_{0}^{\delta} L_1|\lesssim \left(\frac{N_1}{N}\right)^{1-} \|v_1\|_{L^2}
\|(v_2)_x\|_{L^4} \|v_3\|_{L^4}\lesssim N^{-1+}\delta^{\frac{5}{6}-}
N_{\max}^{0-}\prod\limits_{j=1}^3 \|v_i\|_{Y^{1,1/2}}.
\end{equation*}
\end{itemize}

Now, we estimate $\int_{0}^{\delta} L_2$. Our task is to prove that 
\begin{equation}\label{e.aL2}
\begin{split}
&\int_0^{\delta}\sum\limits_{n_1 + n_2 + n_3 = 0} 
\Big|\frac{m(n_1+n_2)-m(n_1)m(n_2)}{m(n_1) m(n_2)}\Big| |n_1+n_2|
\widehat{u_1}(n_1,t) \widehat{u_2}(n_2,t) \widehat{v_3}(n_3,t) \lesssim \\ 
&N^{-1+}\delta^{\frac{19}{24}-}\|u_1\|_{X^{1,1/2}}\|u_2\|_{X^{1,1/2}}
\|v_3\|_{Y^{1,1/2}}
\end{split}
\end{equation}

\begin{itemize}
\item $|n_2|\ll |n_1|\sim |n_3|\gtrsim N$. We estimate the multiplier by 
\begin{equation*}
\Big|\frac{m(n_1+n_2)-m(n_1)m(n_2)}{m(n_1) m(n_2)}\Big|\lesssim
\la(\frac{N_2}{N})^{\frac{1}{2}}\ra.
\end{equation*}
Thus, using $L^2_{xt}L^4_{xt}L^4_{xt}$ H\"older inequality and the
lemmas~\ref{l.Strichartz} and~\ref{l.time-localization}
\begin{equation*}
\begin{split}
\int_0^{\delta}L_2 &\lesssim \la\left(\frac{N_2}{N}\right)^{\frac{1}{2}}\ra
\frac{1}{\la N_2\ra N_3}\delta^{\frac{19}{24}-}
\|u_1\|_{X^{1,1/2}}\|u_2\|_{X^{1,1/2}}
\|v_3\|_{Y^{1,1/2}} \\ 
&\lesssim N^{-1+}\delta^{\frac{19}{24}-} N_{\max}^{0-}
\|u_1\|_{X^{1,1/2}}\|u_2\|_{X^{1,1/2}}\|v_3\|_{Y^{1,1/2}}.
\end{split}
\end{equation*}
\item $|n_1|\ll |n_2|\sim |n_3|$. This case is similar to the previous one.
\item $|n_1|\sim |n_2|\gtrsim N$. Estimating the multiplier by 
\begin{equation*}
\Big|\frac{m(n_1+n_2)-m(n_1)m(n_2)}{m(n_1) m(n_2)}\Big|\lesssim 
\left(\frac{N_2}{N}\right)^{1-}
\end{equation*}
we conclude 
\begin{equation*}
\begin{split}
\int_0^{\delta} L_2&\lesssim \left(\frac{N_2}{N}\right)^{1-} \frac{1}{N_1 N_2}
\delta^{\frac{19}{24}-}\|u_1\|_{X^{1,1/2}}\|u_2\|_{X^{1,1/2}}\|v_3\|_{Y^{1,1/2}}
\\ 
&\lesssim N^{-2+}\delta^{\frac{19}{24}-}N_{\max}^{0-}
\|u_1\|_{X^{1,1/2}}\|u_2\|_{X^{1,1/2}}\|v_3\|_{Y^{1,1/2}}.
\end{split}
\end{equation*}
\end{itemize}

Next, let us compute $\int_0^{\delta} L_3$. We claim that 
\begin{equation}\label{e.aL3}
\begin{split}
&\int_0^{\delta}\sum\limits_{n_1+n_2+n_3=0} 
\Big|\frac{m(n_1+n_2) - m(n_1) m(n_2)}{m(n_1) m(n_2)}\Big| \widehat{u_1}(n_1,t)
\widehat{v_2}(n_2,t) |n_3|\widehat{u_3}(n_3,t) \\ 
&\lesssim N^{-2+}\delta^{\frac{19}{24}-}
\|u_1\|_{X^{1,1/2}}\|v_2\|_{Y^{1,1/2}}\|u_3\|_{X^{1,1/2}}
\end{split}
\end{equation}

\begin{itemize}
\item $|n_2|\ll |n_1|\sim |n_3|$, $|n_1|\gtrsim N$. The multiplier is bounded by
\begin{equation*}
\begin{cases}
\big|\frac{m(n_1+n_2) - m(n_1) m(n_2)}{m(n_1) m(n_2)}\big|\lesssim 
|\frac{\nabla m(n_1)\cdot n_2}{m(n_1)}|\lesssim \frac{N_2}{N_1}, 
\textrm{ if } |n_2|\leq N, \text{ and }\\
\big|\frac{m(n_1+n_2) - m(n_1) m(n_2)}{m(n_1) m(n_2)}\big|\lesssim 
\left(\frac{N_2}{N}\right)^{1/2}, \textrm{ if } |n_2|\geq N.
\end{cases} 
\end{equation*}
So, it is not hard to see that 
\begin{equation*}
\int_0^{\delta} L_3 \lesssim N^{-2+}\delta^{\frac{19}{24}-} N_{\max}^{0-} 
\|u_1\|_{X^{1,1/2}}\|v_2\|_{Y^{1,1/2}}\|u_3\|_{X^{1,1/2}}
\end{equation*}
\item $|n_1|\ll |n_2|\sim |n_3|$, $|n_2|\gtrsim N$. This case is completely 
similar to the previous one. 
\item $|n_1|\sim |n_2|\gtrsim N$. Since the multiplier is bounded by $N_2 / N$,
we get 
\begin{equation*}
\int_0^{\delta} L_3 \lesssim 
N^{-2+}\delta^{\frac{19}{24}-} N_{\max}^{0-} 
\|u_1\|_{X^{1,1/2}}\|v_2\|_{Y^{1,1/2}}\|u_3\|_{X^{1,1/2}}.
\end{equation*}
\end{itemize}

Finally, it remains to estimate the contribution of $\int_0^{\delta} L_4$. It
suffices to see that 
\begin{equation}\label{e.aL4}
\begin{split}
&\int_0^{\delta}\sum\limits_{n_1 + n_2 + n_3 + n_4= 0} 
\Big|\frac{m(n_1+n_2+n_3)-m(n_1)m(n_2)m(n_3)}{m(n_1) m(n_2) m(n_3)}
\Big| |n_4|
\prod\limits_{j=1}^4 \widehat{u_j}(n_j,t) \lesssim \\ 
&N^{-2+}\delta^{\frac{1}{2}-} 
\prod\limits_{j=1}^4\|u_j\|_{X^{1,1/2}}
\end{split}
\end{equation}

\begin{itemize}
\item $N_1, N_2, N_3\gtrsim N$. Since the multiplier verifies 
\begin{equation*}
\Big|\frac{m(n_1+n_2+n_3)-m(n_1)m(n_2)m(n_3)}{m(n_1) m(n_2) m(n_3)}\Big|
\lesssim \left(\frac{N_1}{N}\frac{N_2}{N}\frac{N_3}{N}\right)^{\frac{1}{2}},
\end{equation*}
the application of $L^4_{xt}L^4_{xt}L^4_{xt}L^4_{xt}$ H\"older inequality and
the lemmas~\ref{l.Strichartz},~\ref{l.time-localization} yields 
\begin{equation*}
\int_0^{\delta} L_4\lesssim
\left(\frac{N_1}{N}\frac{N_2}{N}\frac{N_3}{N}\right)^{\frac{1}{2}}
\frac{\delta^{\frac{1}{2}-}}{N_1 N_2 N_3} \prod\limits_{j=1}^4
\|u_j\|_{X^{1,1/2}}\lesssim N^{-3+}\delta^{\frac{1}{2}-} N_{\max}^{0-} 
\prod\limits_{j=1}^4 \|u_j\|_{X^{1,1/2}}.
\end{equation*}
\item $N_1\sim N_2\gtrsim N$ and $N_3, N_4\ll N_1, N_2$. Here the multiplier is
bounded by $\left(\frac{N_1}{N}\frac{N_2}{N}\right)^{\frac{1}{2}} \la
\left(\frac{N_3}{N}\right)^{\frac{1}{2}}\ra$. Hence, 
\begin{equation*}
\int_0^{\delta} L_4\lesssim \left(\frac{N_1}{N}\frac{N_2}{N}\right)^{\frac{1}{2}} \la
\left(\frac{N_3}{N}\right)^{\frac{1}{2}}\ra 
\frac{\delta^{\frac{1}{2}-}}{N_1 N_2 \la N_3\ra} 
\prod\limits_{j=1}^4 \|u_j\|_{X^{1,1/2}}\lesssim 
N^{-2+}\delta^{\frac{1}{2}-} N_{\max}^{0-} 
\prod\limits_{j=1}^4 \|u_j\|_{X^{1,1/2}}.
\end{equation*}
\item $N_1\sim N_4\gtrsim N$ and $N_2, N_3\ll N_1, N_4$. In this case we have
the following estimates for the multiplier 
\begin{equation*}
\Big|\frac{m(n_1+n_2+n_3)-m(n_1)m(n_2)m(n_3)}{m(n_1) m(n_2) m(n_3)}
\Big| \lesssim 
\begin{cases}
\big|\frac{\nabla m(n_1) (n_2+n_3)}{m(n_1)}\big|\lesssim \frac{N_2+N_3}{N_1},
\textrm{ if } 
N_2, N_3\leq N\\
\left(\frac{N_1}{N}\frac{N_2}{N}\right)^{\frac{1}{2}}
\la (\frac{N_3}{N})^{\frac{1}{2}} \ra, \textrm{ if } N_2\geq N,\\
\left(\frac{N_1}{N}\frac{N_3}{N}\right)^{\frac{1}{2}}
\la (\frac{N_2}{N})^{\frac{1}{2}} \ra, \textrm{ if } N_3\geq N.
\end{cases} 
\end{equation*}
Therefore, it is not hard to see that, in any of the situations $N_2, N_3\leq N$,
$N_2\geq N$ or $N_3\geq N$, we have 
\begin{equation*}
\int_0^{\delta} L_4\lesssim N^{-2+}\delta^{\frac{1}{2}-} N_{\max}^{0-}
\prod\limits_{j=1}^4 \|u_j\|_{X^{1,1/2}}.
\end{equation*}
\item $N_1\sim N_2\sim N_4\gtrsim N$ and $N_3\ll N_1, N_2, N_4$. Here we have
the following bound 
\begin{equation*}
\int_0^{\delta} L_4 \lesssim 
\left(\frac{N_1}{N}\frac{N_2}{N}\right)^{\frac{1}{2}}
\la\left(\frac{N_3}{N}\right)^{\frac{1}{2}}\ra \frac{\delta^{\frac{1}{2}-}}{N_1
N_2 N_3} \prod\limits_{j=1}^4 \|u_j\|_{X^{1,1/2}}.
\end{equation*}
\end{itemize}

At this point, clearly the bounds~(\ref{e.aL1}),~(\ref{e.aL2}),~(\ref{e.aL3}) 
and~(\ref{e.aL4}) concludes the proof of the proposition~\ref{p.al}.
\end{proof}

\subsection{Estimates for the modified E-functional}\label{s.e}

\begin{proposition}\label{p.ae}Let $(u,v)$ be a solution of~(\ref{e.nls-kdv}) on
the time interval $[0,\delta]$ such that $\int_{\mathbb{T}}v =0$. Then, for any
$N\geq 1$, $s>1/2$, 
\begin{equation}\label{e.aE}
\begin{split}
&|E(Iu(\delta), Iv(\delta))-E(Iu(0), Iv(0))|\lesssim \\
&\left(N^{-1+}\delta^{\frac{1}{6}-}+N^{-\frac{2}{3}+}\delta^{\frac{3}{8}-} +  
N^{-\frac{3}{2}+}\delta^{\frac{1}{8}-} \right) (\|Iu\|_{X^1}+ \|Iv\|_{Y^1})^3 
+ \\ 
&N^{-1+}\delta^{\frac{1}{2}-} (\|Iu\|_{X^1}+\|Iv\|_{Y^1})^4 + 
N^{-2+}\delta^{\frac{1}{2}-}\|Iu\|_{X^1}^4 (\|Iu\|_{X^1}^2+\|Iv\|_{Y^1}).
\end{split}
\end{equation}
\end{proposition}

\begin{proof}Again we integrate~(\ref{e.ae}) with respect to $t\in [0,\delta]$,
decompose the frequencies into dyadic blocks, etc., so that our objective is to
bound the (integral over $[0,\delta]$ of the) $E_j$ for each $j=1,\dots, 12$. 

For the expression $\int_0^{\delta} E_1$, apply the lemma~\ref{l.duality}. We 
obtain 
\begin{equation*}
|\int_0^{\delta} E_1|\lesssim \|Iv_{xx}\|_{Y^{-1}} \|Iv Iv_x - I(v v_x)\|_{W^1}
\lesssim \|Iv\|_{Y^1} \|Iv Iv_x - I(v v_x)\|_{W^1}
\end{equation*} 
Writing the definition of the norm $W^1$, it suffices to prove the bound  
\begin{equation}
\begin{split}
&\left\|\frac{\la n_3 \ra}{\la\tau_3- n_3^3\ra^{\frac{1}{2}}}
\int \sum \frac{m(n_1+n_2) -
m(n_1)m(n_2)}{m(n_1)m(n_2)} \widehat{v_1}(n_1,\tau_1) \ n_2 \ 
\widehat{v_2}(n_2,\tau_2)\right\|_{L^2_{n_3,\tau_3}} + \\ 
&\left\|\frac{\la n_3 \ra}{\la\tau_3- n_3^3\ra}
\int \sum \frac{m(n_1+n_2) -
m(n_1)m(n_2)}{m(n_1)m(n_2)} \widehat{v_1}(n_1,\tau_1) \ n_2 \ 
\widehat{v_2}(n_2,\tau_2)\right\|_{L^2_{n_3}L^1_{\tau_3}} \lesssim \\ 
& N^{-1+}\delta^{\frac{1}{6}-}\|v_1\|_{Y^{1,1/2}}\|v_2\|_{Y^{1,1/2}}.
\end{split}
\end{equation}
Recall that the dispersion relation $\sum\limits_{j=1}^3 \tau_j
-n_j^3 = -3n_1 n_2 n_3$ implies that, since $n_1 n_2 n_3\neq
0$, if we put $L_j:=|\tau_j - n_j^3|$ and $L_{\max} = \max\{L_j;
j=1,2,3\}$, then $L_{\max}\gtrsim \la n_1\ra \la n_2\ra \la n_3\ra$. 
\begin{itemize}
\item $|n_2|\sim |n_3|\gtrsim N$, $|n_1|\ll |n_2|$. The multiplier is bounded by
\begin{equation*}
\Big|\frac{m(n_1+n_2) - m(n_1)m(n_2)}{m(n_1)m(n_2)}\Big|\lesssim 
\begin{cases}
\frac{N_1}{N_2}, \textrm{ if } |n_1|\leq N, \\ 
\left(\frac{N_1}{N}\right)^{\frac{1}{2}}, \textrm{ if } |n_1|\geq N.
\end{cases}
\end{equation*}
Thus, if $|\tau_3- n_3^3| = L_{\max}$, we have 
\begin{equation*}
\begin{split}
&\left\|\frac{\la n_3 \ra}{\la\tau_3- n_3^3\ra^{\frac{1}{2}}}
\int \sum \frac{m(n_1+n_2) -
m(n_1)m(n_2)}{m(n_1)m(n_2)} \widehat{v_1}(n_1,\tau_1) \ n_2 \ 
\widehat{v_2}(n_2,\tau_2)\right\|_{L^2_{n_3,\tau_3}} \\ 
&\lesssim 
\begin{cases}
\frac{N_1}{N_2}\frac{N_3}{(N_1 N_2 N_3)^{\frac{1}{2}}}
\|v_1\|_{L_{xt}^4}\|(v_2)_x\|_{L^4_{xt}}
\lesssim N^{-1+}
\delta^{\frac{1}{3}-} N_{\max}^{0-} \|v_1\|_{Y^{1,1/2}} \|v_2\|_{Y^{1,1/2}},
\textrm{ if } |n_1|\leq N, \\ 
\left(\frac{N_1}{N_2}\right)^{\frac{1}{2}}\frac{N_3}{N_1}
\frac{1}{(N_1 N_2 N_3)^{\frac{1}{2}}}
\|v_1\|_{L_{xt}^4}\|(v_2)_x\|_{L^4_{xt}}
\lesssim N^{-\frac{3}{2}+}
\delta^{\frac{1}{3}-} N_{\max}^{0-} \|v_1\|_{Y^{1,\frac{1}{2}}} 
\|v_2\|_{Y^{1,\frac{1}{2}}}, 
\textrm{ if } |n_1|\geq N.
\end{cases}
\end{split}
\end{equation*}
and 
\begin{equation*}
\begin{split}
&\left\|\frac{\la n_3 \ra}{\la\tau_3- n_3^3\ra}
\int \sum \frac{m(n_1+n_2) -
m(n_1)m(n_2)}{m(n_1)m(n_2)} \widehat{v_1}(n_1,\tau_1) \ n_2 \ 
\widehat{v_2}(n_2,\tau_2)\right\|_{L^2_{n_3}L^1_{\tau_3}} \\ 
&\lesssim 
\begin{cases}
\frac{N_1}{N_2}\frac{N_3}{(N_1 N_2 N_3)^{\frac{1}{2}-}}\|v_1\|_{L_{xt}^4}
\|(v_2)_x\|_{L^4_{xt}}\lesssim N^{-1+}\delta^{\frac{1}{3}-}N_{\max}^{0-} 
\|v_1\|_{Y^{1,1/2}} \|v_2\|_{Y^{1,1/2}}, 
\textrm{ if } |n_1|\leq N, \\ 
\left(\frac{N_1}{N_2}\right)^{\frac{1}{2}} \frac{N_3}{N_1} 
\frac{\delta^{\frac{1}{3}-}}{(N_1 N_2
N_3)^{\frac{1}{2}-}} \|v_1\|_{Y^{1,\frac{1}{2}}} \|v_2\|_{Y^{1,\frac{1}{2}}}
\lesssim N^{-\frac{3}{2}+}\delta^{\frac{1}{3}-}N_{\max}^{0-} 
\|v_1\|_{Y^{1,\frac{1}{2}}} \|v_2\|_{Y^{1,\frac{1}{2}}}, 
\textrm{ if } |n_1|\geq N.
\end{cases}
\end{split}
\end{equation*}
If either $|\tau_1-n_1^3|=L_{\max}$ or $|\tau_2-n_2^3|=L_{\max}$, we have
\begin{equation*}
\begin{split} 
&\left\|\frac{\la n_3 \ra}{\la\tau_3- n_3^3\ra^{\frac{1}{2}}}
\int \sum \frac{m(n_1+n_2) -
m(n_1)m(n_2)}{m(n_1)m(n_2)} \widehat{v_1}(n_1,\tau_1) \ n_2 \ 
\widehat{v_2}(n_2,\tau_2)\right\|_{L^2_{n_3,\tau_3}} \\ 
&\lesssim 
\begin{cases}
\frac{N_1}{N_2}\frac{N_3}{(N_1 N_2 N_3)^{\frac{1}{2}}}
\frac{\delta^{\frac{1}{6}-}}{N_1}
\|v_1\|_{Y^{1,\frac{1}{2}}}\|v_2\|_{Y^{1,\frac{1}{2}}}
\lesssim N^{-1+}
\delta^{\frac{1}{6}-} N_{\max}^{0-} \|v_1\|_{Y^{1,1/2}} \|v_2\|_{Y^{1,1/2}},
\textrm{ if } |n_1|\leq N, \\ 
\left(\frac{N_1}{N_2}\right)^{\frac{1}{2}}\frac{N_3}{N_1}
\frac{1}{(N_1 N_2 N_3)^{\frac{1}{2}}}
\delta^{\frac{1}{6}-}\|v_1\|_{Y^{1,\frac{1}{2}}} 
\|v_2\|_{Y^{1,\frac{1}{2}}}
\lesssim N^{-\frac{3}{2}+}
\delta^{\frac{1}{3}-} N_{\max}^{0-} \|v_1\|_{Y^{1,\frac{1}{2}}} 
\|v_2\|_{Y^{1,\frac{1}{2}}}, 
\textrm{ if } |n_1|\geq N.
\end{cases}
\end{split}
\end{equation*}
and 
\begin{equation*}
\begin{split}
&\left\|\frac{\la n_3 \ra}{\la\tau_3- n_3^3\ra}
\int \sum \frac{m(n_1+n_2) -
m(n_1)m(n_2)}{m(n_1)m(n_2)} \widehat{v_1}(n_1,\tau_1) \ n_2 \ 
\widehat{v_2}(n_2,\tau_2)\right\|_{L^2_{n_3}L^1_{\tau_3}} \\ 
&\lesssim 
\begin{cases}
\frac{N_1}{N_2}\frac{N_3}{(N_1 N_2
N_3)^{\frac{1}{2}-}}\frac{\delta^{\frac{1}{6}-}}{N_1}
\|v_1\|_{Y^{1,\frac{1}{2}}} \|v_2\|_{Y^{1,\frac{1}{2}}}
\lesssim N^{-1+}\delta^{\frac{1}{6}-}N_{\max}^{0-} 
\|v_1\|_{Y^{1,1/2}} \|v_2\|_{Y^{1,1/2}}, 
\textrm{ if } |n_1|\leq N, \\ 
\left(\frac{N_1}{N_2}\right)^{\frac{1}{2}} \frac{N_3}{N_1}
\frac{\delta^{\frac{1}{6}-}}{(N_1 N_2
N_3)^{\frac{1}{2}-}} \|v_1\|_{Y^{1,\frac{1}{2}}} \|v_2\|_{Y^{1,\frac{1}{2}}}
\lesssim N^{-\frac{3}{2}+}\delta^{\frac{1}{6}-}N_{\max}^{0-} 
\|v_1\|_{Y^{1,\frac{1}{2}}} \|v_2\|_{Y^{1,\frac{1}{2}}}, 
\textrm{ if } |n_1|\geq N.
\end{cases}
\end{split}
\end{equation*}
\item $|n_1|\sim |n_2|\gtrsim N$. Estimating the multiplier by
\begin{equation*}
\Big|\frac{m(n_1+n_2) -
m(n_1)m(n_2)}{m(n_1)m(n_2)}\Big|\lesssim \left(\frac{N_1}{N}\right)^{1-}, 
\end{equation*} 
we have that, if $|\tau_3-n_3^3|=L_{\max}$, 

\begin{equation*}
\begin{split}
&\left\|\frac{\la n_3 \ra}{\la\tau_3- n_3^3\ra^{\frac{1}{2}}}
\int \sum \frac{m(n_1+n_2) -
m(n_1)m(n_2)}{m(n_1)m(n_2)} \widehat{v_1}(n_1,\tau_1) \ n_2 \ 
\widehat{v_2}(n_2,\tau_2)\right\|_{L^2_{n_3,\tau_3}} + \\ 
&\left\|\frac{\la n_3 \ra}{\la\tau_3- n_3^3\ra}
\int \sum \frac{m(n_1+n_2) -
m(n_1)m(n_2)}{m(n_1)m(n_2)} \widehat{v_1}(n_1,\tau_1) \ n_2 \ 
\widehat{v_2}(n_2,\tau_2)\right\|_{L^2_{n_3}L^1_{\tau_3}} \\ 
&\lesssim \left\{\left(\frac{N_1}{N}\right)^{1-}
\frac{N_3}{(N_1 N_2 N_3)^{\frac{1}{2}}}
\frac{\delta^{\frac{1}{3}-}}{N_1} +  
\left(\frac{N_1}{N}\right)^{1-}\frac{N_3}{(N_1 N_2 N_3)^{\frac{1}{2}-}}
\frac{\delta^{\frac{1}{3}-}}{N_1}\right\}
\|v_1\|_{Y^{1,\frac{1}{2}}} \|v_2\|_{Y^{1,\frac{1}{2}}} \\ 
&\lesssim N^{-\frac{3}{2}+}\delta^{\frac{1}{3}-}N_{\max}^{0-}
\|v_1\|_{Y^{1,\frac{1}{2}}} \|v_2\|_{Y^{1,\frac{1}{2}}}
\end{split}
\end{equation*}
and, if either $|\tau_1-n_1^3|=L_{\max}$ or $|\tau_2-n_2^3|=L_{\max}$,
\begin{equation*}
\begin{split}
&\left\|\frac{\la n_3 \ra}{\la\tau_3- n_3^3\ra^{\frac{1}{2}}}
\int \sum \frac{m(n_1+n_2) -
m(n_1)m(n_2)}{m(n_1)m(n_2)} \widehat{v_1}(n_1,\tau_1) \ n_2 \ 
\widehat{v_2}(n_2,\tau_2)\right\|_{L^2_{n_3,\tau_3}} + \\ 
&\left\|\frac{\la n_3 \ra}{\la\tau_3- n_3^3\ra}
\int \sum \frac{m(n_1+n_2) -
m(n_1)m(n_2)}{m(n_1)m(n_2)} \widehat{v_1}(n_1,\tau_1) \ n_2 \ 
\widehat{v_2}(n_2,\tau_2)\right\|_{L^2_{n_3}L^1_{\tau_3}} \\ 
&\lesssim \left\{\left(\frac{N_1}{N}\right)^{1-}
\frac{N_3}{(N_1 N_2 N_3)^{\frac{1}{2}}}\frac{\delta^{\frac{1}{6}-}}{N_1}
+ \left(\frac{N_1}{N}\right)^{1-}
\frac{N_3}{(N_1 N_2 N_3)^{\frac{1}{2}-}}\frac{\delta^{\frac{1}{6}-}}{N_1}
\right\}\|v_1\|_{Y^{1,\frac{1}{2}}} \|v_2\|_{Y^{1,\frac{1}{2}}} \\ 
&\lesssim N^{-\frac{3}{2}+}\delta^{\frac{1}{6}-}N_{\max}^{0-}
\|v_1\|_{Y^{1,\frac{1}{2}}} \|v_2\|_{Y^{1,\frac{1}{2}}}.
\end{split}
\end{equation*}
\end{itemize}

For the expression $\int_0^{\delta}E_2$, it suffices to prove that 
\begin{equation}\label{e.aE2}
\begin{split}
&|\int_0^{\delta}\sum \frac{m(n_3+n_4)-m(n_3)m(n_4)}{m(n_3)m(n_4)}
\widehat{v_1}(n_1,t)\widehat{v_2}(n_2,t)\widehat{v_3}(n_3,t) \ n_4 \ 
\widehat{v_4}(n_4,t)| \lesssim \\ 
& N^{-2+}\delta^{\frac{2}{3}-}\prod\limits_{j=1}^4\|v_j\|_{Y^{1,1/2}}.
\end{split}
\end{equation} 
Since at least two of the $N_i$ are $\gtrsim N$, we can assume that $N_1\geq
N_2\geq N_3$ and $N_1\gtrsim N$. Hence,  
\begin{equation*}
\begin{split}
&\int_0^{\delta} E_2\lesssim \\ 
&\begin{cases}
\left(\frac{N_1}{N}\right)^{1-}\frac{\delta^{\frac{2}{3}-}}{N_1 N_2 N_3}
\prod\limits_{j=1}^4 \|v_j\|_{Y^{1,1/2}}\lesssim N^{-2+}\delta^{\frac{2}{3}-}
N_{\max}^{0-}\prod\limits_{j=1}^4 \|v_j\|_{Y^{1,1/2}}, \textrm{ if } 
|n_3|\sim |n_4|\gtrsim N, \\
\frac{N_3}{N_4}\frac{\delta^{\frac{2}{3}-}}{N_1 N_2 N_3}
\prod\limits_{j=1}^4 \|v_j\|_{Y^{1,1/2}}\lesssim N^{-2+}\delta^{\frac{2}{3}-}
N_{\max}^{0-}\prod\limits_{j=1}^4 \|v_j\|_{Y^{1,1/2}}, \textrm{ if } 
|n_3|\ll |n_4|, |n_3|\leq N |n_4|\gtrsim N, \\
\left(\frac{N_3}{N}\right)^{\frac{1}{2}}
\frac{\delta^{\frac{2}{3}-}}{N_1 N_2 N_3}
\prod\limits_{j=1}^4 \|v_j\|_{Y^{1,\frac{1}{2}}}
\lesssim N^{-2+}\delta^{\frac{2}{3}-}
N_{\max}^{0-}\prod\limits_{j=1}^4 \|v_j\|_{Y^{1,\frac{1}{2}}}, \textrm{ if } 
|n_3|\ll |n_4|, |n_3|\geq N, |n_4|\gtrsim N. \\
\end{cases}
\end{split}
\end{equation*}

Next, we estimate the contribution of $\int_0^{\delta} E_3$. We claim that 
\begin{equation}\label{l.aE3}
\begin{split}
&\int_{0}^{\delta}\sum\frac{m(n_1n_2n_3) - 
m(n_1)m(n_2)m(n_3)}{m(n_1)m(n_2)m(n_3)} \widehat{u_1}(n_1,t)
\widehat{u_2}(n_2,t) \widehat{u_3}(n_3,t)\ |n_4|^2 \
\widehat{u_4}(n_4,t)\lesssim \\ 
& N^{-1+}\delta^{\frac{1}{2}-}\prod\limits_{j=1}^4 \|u_j\|_{X^{1,1/2}}.
\end{split}
\end{equation}
\begin{itemize}
\item $|n_1|\sim |n_2|\sim |n_3|\sim |n_4|\gtrsim N$. Since the multiplier
satisfies 
\begin{equation*}
\frac{m(n_1n_2n_3) - 
m(n_1)m(n_2)m(n_3)}{m(n_1)m(n_2)m(n_3)}\lesssim
\left(\frac{N_1}{N}\right)^{\frac{3}{2}}
\end{equation*}
we obtain 
\begin{equation*}
\int_0^{\delta} E_3\lesssim \left(\frac{N_1}{N}\right)^{\frac{3}{2}}
\frac{N_4}{N_1N_2N_3}\delta^{\frac{1}{2}-} 
\prod\limits_{j=1}^4 \|u_j\|_{X^{1,1/2}}\lesssim N^{-2+}\delta^{\frac{1}{2}-}
N_{\max}^{0-}\prod\limits_{j=1}^4 \|u_j\|_{X^{1,1/2}}.
\end{equation*}
\item Exactly two frequencies are $\gtrsim N$. We consider the most difficult
case $|n_4|\gtrsim N$, $|n_1|\sim |n_4|$ and $|n_2|, |n_3|\ll |n_1|, |n_4|$. The
multiplier is estimated by 
\begin{equation*}
\frac{m(n_1n_2n_3) - 
m(n_1)m(n_2)m(n_3)}{m(n_1)m(n_2)m(n_3)}\lesssim
\begin{cases}
\la\left(\frac{N_3}{N}\right)^{\frac{1}{2}}\ra 
\left(\frac{N_2}{N}\right)^{\frac{1}{2}}, \textrm{ if } |n_2|\geq N,\\ 
\la\left(\frac{N_2}{N}\right)^{\frac{1}{2}}\ra 
\left(\frac{N_3}{N}\right)^{\frac{1}{2}}, \textrm{ if } |n_3|\geq N,\\ 
\frac{N_2+N_3}{N_1}, \textrm{ if } |n_2|, |n_3|\leq N.
\end{cases}
\end{equation*}
Thus, 
\begin{equation*}
\int_0^{\delta}E_3\lesssim N^{-1+}\delta^{\frac{1}{2}-}
N_{\max}^{0-}\prod\limits_{j=1}^4 \|u_j\|_{X^{1,1/2}}.
\end{equation*}
\item Exactly three frequencies are $\gtrsim N$. The most difficult case is
$|n_1|\sim |n_2|\sim |n_4|\gtrsim N$ and $|n_3|\ll |n_1|, |n_2|, |n_4|$. Here
the multiplier is bounded by 
\begin{equation*}
\frac{m(n_1n_2n_3) - 
m(n_1)m(n_2)m(n_3)}{m(n_1)m(n_2)m(n_3)}\lesssim 
\left(\frac{N_1}{N}\frac{N_2}{N}\right)^{\frac{1}{2}}
\la\left(\frac{N_3}{N}\right)^{\frac{1}{2}}\ra.
\end{equation*}
Hence, 
\begin{equation*}
\int_0^{\delta} E_3\lesssim 
\left(\frac{N_1}{N}\frac{N_2}{N}\right)^{\frac{1}{2}}
\la\left(\frac{N_3}{N}\right)^{\frac{1}{2}}\ra \frac{N_4}{N_1 N_2 N_3} 
\delta^{\frac{1}{2}-}\prod\limits_{j=1}^4 \|u_j\|_{X^{1,1/2}}\lesssim 
N^{-1+}\delta^{\frac{1}{2}-}
N_{\max}^{0-}\prod\limits_{j=1}^4 \|u_j\|_{X^{1,1/2}}.
\end{equation*}
\end{itemize}

The contribution of $\int_0^{\delta} E_4$ is controlled if we are able to show
that 
\begin{equation}\label{e.aE4}
\begin{split}
&\int_0^{\delta}\sum \frac{m(n_1+n_2) - m(n_1)m(n_2)}{m(n_1)m(n_2)} 
\widehat{v_1}(n_1,t)\ |n_2| \ \widehat{v_2}(n_2,t) \widehat{u_3}(n_3,t)
\widehat{u_4}(n_4,t)\lesssim \\ 
& N^{-1+}\delta^{\frac{7}{12}-}\prod\limits_{j=1}^2
\|u_j\|_{X^{1,1/2}}\|v_j\|_{Y^{1,1/2}}.
\end{split}
\end{equation}
We crudely bound the multiplier by 
\begin{equation*}
|\frac{m(n_1+n_2) - m(n_1)m(n_2)}{m(n_1)m(n_2)}|\lesssim
\left(\frac{N_{\max}}{N}\right)^{1-}.
\end{equation*}
The most difficult case is $|n_2|\geq N$. We have two possibilities: 
\begin{itemize}
\item Exactly two frequencies are $\gtrsim N$. We can assume $N_3\ll N_2$. In
particular, 
\begin{equation*}
\int_0^{\delta}E_4\lesssim \left(\frac{N_{\max}}{N}\right)^{1-}
\frac{\delta^{\frac{7}{12}}}{N_1 N_3 N_4}\prod\limits_{j=1}^2
\|u_j\|_{X^{1,\frac{1}{2}}}\|v_j\|_{Y^{1,\frac{1}{2}}}\lesssim 
N^{-1+}\delta^{\frac{7}{12}-}N_{\max}^{0-}\prod\limits_{j=1}^2
\|u_j\|_{X^{1,\frac{1}{2}}}\|v_j\|_{Y^{1,\frac{1}{2}}}.
\end{equation*}
\item At least three frequencies are $\gtrsim N$. In this case, 
\begin{equation*}
\int_0^{\delta}E_4\lesssim 
N^{-2+}\delta^{\frac{7}{12}-}N_{\max}^{0-}\prod\limits_{j=1}^2
\|u_j\|_{X^{1,\frac{1}{2}}}\|v_j\|_{Y^{1,\frac{1}{2}}}.
\end{equation*}
\end{itemize}

The expression $\int_0^{\delta}E_5$ is controlled if we are able to prove 
\begin{equation}\label{e.aE5}
\begin{split}
&\int_0^{\delta}\sum \frac{m(n_1+n_2) - m(n_1)m(n_2)}{m(n_1)m(n_2)} 
\widehat{u_1}(n_1,t) \widehat{u_2}(n_2,t) \widehat{v_3}(n_3,t)
\ |n_4| \ \widehat{v_4}(n_4,t)\lesssim \\ 
& N^{-1+}\delta^{\frac{7}{12}-}\prod\limits_{j=1}^2
\|u_j\|_{X^{1,1/2}}\|v_j\|_{Y^{1,1/2}}.
\end{split}
\end{equation}
This follows directly from the previous analysis for~(\ref{e.aE4}).

For the term $\int_0^{\delta} E_6$, we apply the lemma~\ref{l.duality} 
to obtain 
\begin{equation*}
\int_{0}^{\delta}E_6\lesssim \|(Iv)_{xx}\|_{Y^{-1}}\|(|Iu|^2 -
I(|u|^2))_x\|_{W^1}\lesssim \|Iv\|_{Y^1}\|(|Iu|^2 -
I(|u|^2))_x\|_{W^1}.
\end{equation*}
So, the definition of the $W^1$ norm means that we have to prove 
\begin{equation}\label{e.aE6}
\begin{split}
&\left\|\frac{\la n_3\ra}{\la\tau_3-n_3^3\ra^{\frac{1}{2}}} |n_3| \int\sum
\frac{m(n_1+n_2) - m(n_1)m(n_2)}{m(n_1)m(n_2)} \widehat{u_1}(n_1,\tau_1)
\widehat{u_2}(n_2,\tau_2)\right\|_{L_{n_3,\tau_3}^2} + \\ 
&\left\|\frac{\la n_3\ra}{\la\tau_3-n_3^3\ra^{\frac{1}{2}}} |n_3| \int\sum
\frac{m(n_1+n_2) - m(n_1)m(n_2)}{m(n_1)m(n_2)} \widehat{u_1}(n_1,\tau_1)
\widehat{u_2}(n_2,\tau_2)\right\|_{L^2_{n_3}L_{\tau_3}^1}\lesssim \\ 
& \left\{N^{-\frac{3}{2}+}\delta^{\frac{1}{8}-} +
N^{-\frac{2}{3}}\delta^{\frac{3}{8}-}\right\}
\|u_1\|_{X^{1,1/2}}\|u_2\|_{X^{1,1/2}}.
\end{split}
\end{equation}
Note that $\sum\tau_j = 0$ and $\sum n_j = 0$. In particular, we obtain the
dispersion relation 
\begin{equation*}
\tau_3 - n_3^3 + \tau_2 + n_2^2 + \tau_1 + n_1^2 = -n_3^3+n_1^2+n_2^2.
\end{equation*}
\begin{itemize}
\item $|n_1|\gtrsim N$, $|n_2|\ll |n_1|$. Denoting by $L_1:=|\tau_1+n_1^2|$,
$L_2:=|\tau_2+n_2^2|$ and $L_3:=|\tau_3-n_3^3|$, the dispersion relation says
that in the present situation $L_{\max}:=\max\{L_j\}\gtrsim N_3^3$. Since the
multiplier is bounded by 
\begin{equation*}
\Big|\frac{m(n_1+n_2) - m(n_1)m(n_2)}{m(n_1)m(n_2)}\Big|\lesssim 
\begin{cases}
\frac{\nabla m(n_1) n_2}{m(n_1)}\lesssim \frac{N_2}{N_1}, \textrm{ if }
|n_2|\leq N, \\
\left(\frac{N_2}{N}\right)^{\frac{1}{2}}, \textrm{ if } |n_2|\geq N,
\end{cases}
\end{equation*}
we deduce that 
\begin{equation*}
\begin{split}
&\left\|\frac{\la n_3\ra}{\la\tau_3-n_3^3\ra^{\frac{1}{2}}} |n_3| \int\sum
\frac{m(n_1+n_2) - m(n_1)m(n_2)}{m(n_1)m(n_2)} \widehat{u_1}(n_1,\tau_1)
\widehat{u_2}(n_2,\tau_2)\right\|_{L_{n_3,\tau_3}^2} + \\ 
&\left\|\frac{\la n_3\ra}{\la\tau_3-n_3^3\ra^{\frac{1}{2}}} |n_3| \int\sum
\frac{m(n_1+n_2) - m(n_1)m(n_2)}{m(n_1)m(n_2)} \widehat{u_1}(n_1,\tau_1)
\widehat{u_2}(n_2,\tau_2)\right\|_{L^2_{n_3}L_{\tau_3}^1}\lesssim \\ 
& \frac{N_3^2}{N_3^{\frac{3}{2}-}}\frac{\delta^{\frac{1}{8}-}}{N N_1} 
\|u_1\|_{X^{1,1/2}}\|u_2\|_{X^{1,1/2}} \lesssim  
N^{-\frac{3}{2}+}\delta^{\frac{1}{8}-} N_{\max}^{0-}
\|u_1\|_{X^{1,1/2}}\|u_2\|_{X^{1,1/2}}.
\end{split}
\end{equation*}
\item $|n_1|\sim |n_2|\gtrsim N$, $|n_3|^3\gg |n_2|^2$. In the present case the
multiplier is bounded by $\left(\frac{N_1}{N}\right)^{1-}$ and the dispersion 
relation says that $L_{\max}\gtrsim N_3^3$. Thus, 
\begin{equation*}
\begin{split}
&\left\|\frac{\la n_3\ra}{\la\tau_3-n_3^3\ra^{\frac{1}{2}}} |n_3| \int\sum
\frac{m(n_1+n_2) - m(n_1)m(n_2)}{m(n_1)m(n_2)} \widehat{u_1}(n_1,\tau_1)
\widehat{u_2}(n_2,\tau_2)\right\|_{L_{n_3,\tau_3}^2} + \\ 
&\left\|\frac{\la n_3\ra}{\la\tau_3-n_3^3\ra^{\frac{1}{2}}} |n_3| \int\sum
\frac{m(n_1+n_2) - m(n_1)m(n_2)}{m(n_1)m(n_2)} \widehat{u_1}(n_1,\tau_1)
\widehat{u_2}(n_2,\tau_2)\right\|_{L^2_{n_3}L_{\tau_3}^1}\lesssim \\ 
& \frac{N_3^2}{N_3^{\frac{3}{2}-}}\left(\frac{N_1}{N}\right)^{1-}
\frac{\delta^{\frac{1}{8}-}}{N_1 N_2} 
\|u_1\|_{X^{1,1/2}}\|u_2\|_{X^{1,1/2}} \lesssim  
N^{-\frac{3}{2}+}\delta^{\frac{1}{8}-} N_{\max}^{0-}
\|u_1\|_{X^{1,1/2}}\|u_2\|_{X^{1,1/2}}.
\end{split}
\end{equation*} 
\item $|n_1|\sim |n_2|\gtrsim N$ and $|n_3|^3\lesssim |n_2|^2$. Here the
dispersion relation does not give useful information about $L_{\max}$. Since the
multiplier is estimated by $\left(\frac{N_2}{N}\right)^{\frac{1}{2}}$, we obtain
the crude bound 
\begin{equation*}
\begin{split}
&\left\|\frac{\la n_3\ra}{\la\tau_3-n_3^3\ra^{\frac{1}{2}}} |n_3| \int\sum
\frac{m(n_1+n_2) - m(n_1)m(n_2)}{m(n_1)m(n_2)} \widehat{u_1}(n_1,\tau_1)
\widehat{u_2}(n_2,\tau_2)\right\|_{L_{n_3,\tau_3}^2} + \\ 
&\left\|\frac{\la n_3\ra}{\la\tau_3-n_3^3\ra^{\frac{1}{2}}} |n_3| \int\sum
\frac{m(n_1+n_2) - m(n_1)m(n_2)}{m(n_1)m(n_2)} \widehat{u_1}(n_1,\tau_1)
\widehat{u_2}(n_2,\tau_2)\right\|_{L^2_{n_3}L_{\tau_3}^1}\lesssim \\ 
& N_3^2\left(\frac{N_2}{N}\right)^{\frac{1}{2}}
\frac{\delta^{\frac{3}{8}-}}{N_1 N_2} 
\|u_1\|_{X^{1,1/2}}\|u_2\|_{X^{1,1/2}} \lesssim  
N^{-\frac{2}{3}+}\delta^{\frac{3}{8}-} N_{\max}^{0-}
\|u_1\|_{X^{1,1/2}}\|u_2\|_{X^{1,1/2}}.
\end{split}
\end{equation*} 
\end{itemize}

Next, the desired bound related to $\int_0^{\delta} E_7$ follows from 
\begin{equation}\label{e.aE7}
\begin{split}
&\int_0^{\delta}\sum\Big|\frac{m(n_1+n_2)-m(n_1)m(n_2)}{m(n_1)m(n_2)}\Big|
|n_1+n_2|\widehat{u_1}(n_1,t)\widehat{v_2}(n_2,t) |n_3|\widehat{u_3}(n_3,t)
\lesssim \\ 
& N^{-1+}\delta^{\frac{19}{24}-}\|u_1\|_{X^{1,1/2}}\|v_2\|_{Y^{1,1/2}}
\|u_3\|_{X^{1,1/2}}
\end{split}
\end{equation}
\begin{itemize}
\item $|n_1|\ll |n_2|\gtrsim N$. The multiplier is $\lesssim (|n_2|/N)^{1/2}$ so
that 
\begin{equation*}
\begin{split}
&\int_0^{\delta}E_7\lesssim \frac{1}{N^{1/2}} \int_0^{\delta}\sum |n_1+n_2| 
\widehat{u_1}(n_1,t) |n_2|^{1/2} \widehat{v_2}(n_2,t) |n_3|\widehat{u_3}(n_3,t)
\lesssim \\ 
& N^{-1}\delta^{\frac{19}{24}-}\|u_1\|_{X^{1,1/2}}\|v_2\|_{Y^{1,1/2}}
\|u_3\|_{X^{1,1/2}}.
\end{split}
\end{equation*}
\item $|n_1|\sim |n_2|\gtrsim N$. The multiplier is $\lesssim |n_2|/N$. Hence, 
\begin{equation*}
\int_0^{\delta}E_7\lesssim N^{-1}\delta^{\frac{19}{24}-}
\|u_1\|_{X^{1,1/2}}\|v_2\|_{Y^{1,1/2}}\|u_3\|_{X^{1,1/2}}.
\end{equation*}
\item $|n_1|\gtrsim N$, $|n_2|\leq N$. The multiplier is again $\lesssim N_2/N$,
so that it can be estimated as above.
\end{itemize}

Now we turn to the term $\int_0^{\delta}E_8$. The objective is to show that 
\begin{equation}\label{e.aE8}
\int_0^{\delta}\Big|\frac{m(n_1+n_2)-m(n_1)m(n_2)}{m(n_1)m(n_2)}\Big|
|n_1+n_2|\prod\limits_{j=1}^4\widehat{u_j}(n_j,t)\lesssim  
N^{-1+}\delta^{\frac{1}{2}-}\prod\limits_{j=1}^4\|u_j\|_{X^{1,1/2}}
\end{equation}
\begin{itemize}
\item At least three frequencies are $\gtrsim N$. We can assume $|n_1|\geq
|n_2|$. The multiplier is bounded by $N_{\max}/N$ so that 
\begin{equation*}
\int_0^{\delta}E_8\lesssim \frac{N_{\max}}{N}
\frac{\delta^{\frac{1}{2}-}}{N_2 N_3 N_4}
\prod\limits_{j=1}^4\|u_j\|_{X^{1,1/2}}\lesssim 
N^{-2+}\delta^{\frac{1}{2}-}N_{\max}^{0-}
\prod\limits_{j=1}^4\|u_j\|_{X^{1,1/2}}.
\end{equation*} 
\item Exactly two frequencies are $\gtrsim N$. Without loss of generality, we
suppose $|n_1|\sim |n_2|\gtrsim N$ and $|n_3|, |n_4|\ll N$. Since the multiplier
satisfies 
\begin{equation*}
\Big|\frac{m(n_1+n_2)-m(n_1)m(n_2)}{m(n_1)m(n_2)}\Big|\lesssim 
\left(\frac{N_{\max}}{N}\right)^{1-},
\end{equation*}
we get the bound 
\begin{equation*}
\int_0^{\delta}E_8\lesssim \left(\frac{N_{\max}}{N}\right)^{1-}
\frac{\delta^{\frac{1}{2}-}}{N_2 N_3 N_4}
\prod\limits_{j=1}^4\|u_j\|_{X^{1,1/2}}\lesssim 
N^{-1+}\delta^{\frac{1}{2}-}N_{\max}^{0-}
\prod\limits_{j=1}^4\|u_j\|_{X^{1,1/2}}.
\end{equation*}
\end{itemize}

The contribution of $\int_0^{\delta}E_9$ is estimated if we prove that 
\begin{equation}\label{e.aE9}
\begin{split}
&\int_0^{\delta}\Big|\frac{m(n_1+n_2)-m(n_1)m(n_2)}{m(n_1)m(n_2)}\Big|
\widehat{u_1}(n_1,t)\widehat{v_2}(n_2,t)\widehat{u_3}(n_3,t)\widehat{v_4}(n_4,t)
\lesssim \\ 
&N^{-2+}\delta^{\frac{7}{12}-}\|u_1\|_{X^{1,1/2}}\|v_2\|_{Y^{1,1/2}}
\|u_3\|_{X^{1,1/2}}\|v_4\|_{Y^{1,1/2}}.
\end{split}
\end{equation}
This follows since at least two frequencies are $\gtrsim N$ and the multiplier 
is always bounded by $(N_{\max}/N)^{1-}$, so that 
\begin{equation*}
\begin{split}
&\int_0^{\delta}E_9\lesssim \left(\frac{N_{\max}}{N}\right)^{1-} \|u_1\|_{L^4}
\|v_2\|_{L^4}\|u_3\|_{L^4}\|v_4\|_{L^4}\lesssim \\ 
& \left(\frac{N_{\max}}{N}\right)^{1-}
\frac{\delta^{\frac{1}{4}+ \frac{1}{3}-}}{N_1 N_2 N_3 N_4} 
\|u_1\|_{X^{1,1/2}}\|v_2\|_{Y^{1,1/2}}
\|u_3\|_{X^{1,1/2}}\|v_4\|_{Y^{1,1/2}}\lesssim \\ 
& N^{-2+}\delta^{\frac{7}{12}-}\|u_1\|_{X^{1,1/2}}\|v_2\|_{Y^{1,1/2}}
\|u_3\|_{X^{1,1/2}}\|v_4\|_{Y^{1,1/2}}.
\end{split}
\end{equation*} 

Now, we treat the term $\int_0^{\delta}E_{10}$. It is sufficient to prove 
\begin{equation}\label{e.aE10}
\begin{split}
&\int_0^{\delta}\sum
\Big|\frac{m(n_4+n_5+n_6)-m(n_4)m(n_5)m(n_6)}{m(n_4)m(n_5)m(n_6)}\Big|
\prod\limits_{j=1}^6\widehat{u_j}(n_j,t)\lesssim \\ 
& N^{-2+}\delta^{\frac{1}{2}-}\prod\limits_{j=1}^6\|u_j\|_{X^1}
\end{split}
\end{equation}
However, this follows easily from the facts that the multiplier is bounded by
$(N_{\max}/N)^{3/2}$, at least two frequencies are $\gtrsim N$, say
$|n_{i_1}|\geq |n_{i_2}|\gtrsim N$, the Strichartz bound $X^{0,3/8}\subset L^4$
and the inclusion\footnote{This inclusion is an easy consequence of Sobolev 
embedding.} $X^{\frac{1}{2}+}\subset L_{xt}^{\infty}$. Indeed, if we combine
these informations, it is not hard to get 
\begin{equation*}
\begin{split}
\int_0^{\delta}E_{10}\lesssim \left(\frac{N_{\max}}{N}\right)^{\frac{3}{2}}
\frac{1}{N_{i_1} N_{i_2} N_{i_3}
N_{i_4}}\delta^{\frac{1}{2}-}\frac{1}{(N_{i_5}N_{i_6})^{1/2-}}
\prod\limits_{j=1}^6\|u_j\|_{X^1}\lesssim 
N^{-2+}\delta^{\frac{1}{2}-}N_{\max}^{0-}\prod\limits_{j=1}^6\|u_j\|_{X^1} 
\end{split}
\end{equation*}

For the expression $\int_0^{\delta}E_{11}$, we use again that the multiplier is
bounded by $(N_{\max}/N)^{3/2}$, at least two frequencies are $\gtrsim N$ (say
$|n_{i_1}|\geq |n_{i_2}|\gtrsim N$), the
Strichartz bounds in lemma~\ref{l.Strichartz} and the inclusions
$X^{\frac{1}{2}+},Y^{\frac{1}{2}+}\subset L^{\infty}_{xt}$ to obtain 
\begin{equation}
\begin{split}
&\int_0^{\delta}\sum 
\Big|\frac{m(n_1+n_2+n_3)-m(n_1)m(n_2)m(n_3)}{m(n_1)m(n_2)m(n_3)}\Big|
\prod\limits_{j=1}^4\widehat{u_j}(n_j,t)\widehat{v_5}(n_5,t)\lesssim \\ 
& \left(\frac{N_{\max}}{N}\right)^{\frac{3}{2}}
\frac{1}{N_{i_1} N_{i_2} N_{i_3} N_{i_4}}
\frac{\delta^{\frac{1}{2}-}}{N_{i_5}^{1/2-}}
\prod\limits_{j=1}^4\|u_j\|_{X^1} \|v_5\|_{Y^1}\lesssim \\
& N^{-2+}\delta^{\frac{1}{2}-}\prod\limits_{j=1}^4\|u_j\|_{X^1} \|v_5\|_{Y^1}.
\end{split}
\end{equation}

The analysis of $\int_0^{\delta}E_{12}$ is similar to the
$\int_0^{\delta}E_{11}$. This completes the proof of the proposition~\ref{p.ae}.
\end{proof} 

\section{Global well-posedness below the energy space}\label{s.global}

In this section we combine the variant local well-posedness result in
proposition~\ref{p.local} with the two almost conservation results in the
propositions~\ref{p.al} and~\ref{p.ae} to prove the theorem~\ref{t.A}.

\begin{remark} Note that the spatial mean $\int_{\T} v(t,x) dx$ is 
preserved during the evolution~(\ref{e.nls-kdv}). 
Thus, we can assume that the initial data $v_0$ has zero-mean, since otherwise 
we make the change $w= v-\int_{\T}v_0 dx$ at the expense of two harmless linear 
terms (namely, $u\int_{\T}v_0 dx$ and $\p_x v \int_{\T}v_0$).  
\end{remark}

The definition of the I-operator implies that the initial data satisfies
$\|Iu_0\|_{H^1}^2 +\|Iv_0\|_{H^1}^2\lesssim N^{2(1-s)}$ and
$\|Iu_0\|_{L^2}^2 +\|Iv_0\|_{L^2}^2\lesssim 1$. By the estimates~(\ref{e.L1})
and~(\ref{e.E3}), we get that $|L(Iu_0, Iv_0)|\lesssim N^{1-s}$ and $|E(Iu_0,
Iv_0)|\lesssim N^{2(1-s)}$.

Also, any bound for $L(Iu,Iv)$ and $E(Iu,Iv)$ of the form 
$|L(Iu,Iv)|\lesssim N^{1-s}$ and $|E(Iu,Iv)|\lesssim N^{2(1-s)}$ implies that 
$\|Iu\|_{L^2}^2\lesssim M$, $\|Iv\|_{L^2}^2\lesssim N^{1-s}$ and 
$\|Iu\|_{H^1}^2+\|Iv\|_{H^1}^2\lesssim N^{2(1-s)}$.

Given a time $T$, if we can uniformly bound the $H^1$-norms of the solution at
times $t=\delta$, $t=2\delta$, etc., the local existence result in
proposition~\ref{p.local} says that the solution can be extended up to any time
interval where such a uniform bound holds. On the other hand, given a time $T$, 
if we can interact $T\delta^{-1}$ times the local existence result, the 
solution exists in the time interval $[0,T]$. So, in view of the 
propositions~\ref{p.al} and~\ref{p.ae}, it suffices to show  

\begin{equation}\label{e.L} 
(N^{-1+}\delta^{\frac{19}{24}-}N^{3(1-s)} +
N^{-2+}\delta^{\frac{1}{2}-}N^{4(1-s)})T\delta^{-1}\lesssim 
N^{1-s}
\end{equation}
and 
\begin{equation}\label{e.E} 
\begin{split}
&\left\{(N^{-1+}\delta^{\frac{1}{6}-} + N^{-\frac{2}{3}+}\delta^{\frac{3}{8}-} +
N^{-\frac{3}{2}+}\delta^{\frac{1}{8}-})N^{3(1-s)} +
N^{-1+}\delta^{\frac{1}{2}-}N^{4(1-s)} + 
N^{-2+}\delta^{\frac{1}{2}-}N^{6(1-s)}\right\}\frac{T}{\delta}\lesssim \\
& N^{2(1-s)}
\end{split}
\end{equation}

At this point, we recall that the proposition~\ref{p.local} says that
$\delta\sim N^{-\frac{16}{3}(1-s)-}$ if $\beta\neq 0$ and $\delta\sim
N^{-8(1-s)-}$ if $\beta=0$. Hence, 

\begin{itemize}
\item $\beta\neq 0$. The condition~(\ref{e.L}) holds for 
\begin{equation*}
-1+\frac{5}{24}\frac{16}{3}(1-s)+3(1-s)< (1-s), \textrm{ i.e. },
s>19/28
\end{equation*}
and
\begin{equation*}
-2+\frac{1}{2}\frac{16}{3}(1-s)+4(1-s)< (1-s), \textrm{ i.e. },
s>11/17;
\end{equation*}
Similarly, the condition~(\ref{e.E}) is satisfied if 
\begin{equation*}
-1+\frac{5}{6}\frac{16}{3}(1-s) + 3(1-s)< 2(1-s), \textrm{ i.e. }, s>40/49;
\end{equation*}
\begin{equation*}
-\frac{2}{3}+\frac{5}{6}\frac{16}{3}(1-s)+3(1-s)< 2(1-s), \textrm{ i.e. },
s>11/13;
\end{equation*}
\begin{equation*}
-\frac{3}{2}+\frac{7}{8}\frac{16}{3}(1-s)+3(1-s)< 2(1-s), \textrm{ i.e. },
s>25/34;
\end{equation*}
\begin{equation*}
-1+\frac{1}{2}\frac{16}{3}(1-s)+4(1-s)< 2(1-s), \textrm{ i.e. },
s>11/14 
\end{equation*}
and 
\begin{equation*}
-2+\frac{1}{2}\frac{16}{3}(1-s)+6(1-s)< 2(1-s), \textrm{ i.e. },
s>7/10.
\end{equation*}

Thus, we conclude that the non-resonant NLS-KdV system is globally well-posed
for any $s>11/13$. 

\item $\beta=0$. The condition~(\ref{e.L}) is fulfilled when 
\begin{equation*}
-1+\frac{5}{24}8(1-s)+3(1-s)< (1-s), \textrm{ i.e. },
s>8/11
\end{equation*}
and
\begin{equation*}
-2+\frac{1}{2}8(1-s)+4(1-s)< (1-s), \textrm{ i.e. },
s>5/7;
\end{equation*}
Similarly, the condition~(\ref{e.E}) is verified for 
\begin{equation*}
-1+\frac{5}{6}8(1-s) + 3(1-s)< 2(1-s), \textrm{ i.e. }, s>20/23;
\end{equation*}
\begin{equation*}
-\frac{2}{3}+\frac{5}{6}8(1-s)+3(1-s)< 2(1-s), \textrm{ i.e. },
s>8/9;
\end{equation*}
\begin{equation*}
-\frac{3}{2}+\frac{7}{8}8(1-s)+3(1-s)< 2(1-s), \textrm{ i.e. },
s>13/16;
\end{equation*}
\begin{equation*}
-1+\frac{1}{2}8(1-s)+4(1-s)< 2(1-s), \textrm{ i.e. },
s>5/6 
\end{equation*}
and 
\begin{equation*}
-2+\frac{1}{2}8(1-s)+6(1-s)< 2(1-s), \textrm{ i.e. },
s>3/4.
\end{equation*}

Hence, we obtain that the resonant NLS-KdV system is globally well-posed
for any $s>8/9$.
\end{itemize}


 
\end{document}